\documentclass[letterpaper, paper,11pt]{AAS}		% for final proceedings (20-page limit)

\usepackage{bm}
\usepackage{amsmath}
\usepackage{subfigure}
\usepackage[colorlinks=true, pdfstartview=FitV, linkcolor=black, citecolor= black, urlcolor= black]{hyperref}
\usepackage{overcite}
\usepackage{footnpag}			      	% make footnote symbols restart on each page

\RequirePackage{color}
\usepackage{soul}
\usepackage[per-mode = symbol]{siunitx}  % For SI measure units. The [...] grants that m/s is visualized instead of m*s^-1, that would be [per-mode = reciprocal]. E.g. \SI{12}{\W\per\nano\m}
\usepackage{comment}

\PaperNumber{19-763}

\makeatletter
\renewcommand{\thesubfigure}{\alph{subfigure}}% (a) -> a
\renewcommand{\@thesubfigure}{(\thesubfigure)\hskip\subfiglabelskip}% a -> a)
\makeatother

\begin{document}

\title{A Convex Optimization Approach for Finite-Thrust Time-Constrained Cooperative Rendezvous}

\author{Boris Benedikter\thanks{PhD Student, Department of Mechanical and Aerospace Engineering, Sapienza University of Rome, Via Eudossiana 18 - 00184, Rome, Italy},  
Alessandro Zavoli\thanks{Research Assistant, Department of Mechanical and Aerospace Engineering, Sapienza University of Rome, Via Eudossiana 18 - 00184, Rome, Italy},
\ and Guido Colasurdo\thanks{Full Professor, Department of Mechanical and Aerospace Engineering, Sapienza University of Rome, Via Eudossiana 18 - 00184, Rome, Italy}
}

\maketitle{} 		

\begin{abstract}
This paper presents a convex approach to the optimization of a cooperative rendezvous, that is, the problem of two distant spacecraft that simultaneously operate to get closer. 
Convex programming guarantees convergence towards the optimal solution in a limited, short, time by using highly efficient numerical algorithms. 
A combination of lossless and successive convexification techniques is adopted to handle the nonconvexities of the original problem. 
Specifically, a convenient change of variables and a constraint relaxation are performed, while a successive linearization of the equations of motion is employed to handle the nonlinear dynamics. 
A filtering technique concerning the recursive update of the reference solution is proposed in order to enhance the algorithm robustness.
Numerical results are presented and compared with those provided by an indirect method.
\end{abstract}

\section{Introduction}
The rendezvous problem is a traditional, well-established research topic in spaceflight mechanics, as it is a basic maneuver in several operative missions, such as interplanetary exploration, on-orbit servicing or docking, and remote sensing.
A broad literature is thus available, but most of it focuses on the non-cooperative rendezvous problem, where an active chaser maneuvers to reach a passive target.
Instead, this paper proposes a convex programming approach to the optimization of a \emph{cooperative} rendezvous, that is, the problem of two distant spacecraft that simultaneously and coordinately operate to get closer.

Both indirect and direct methods have been proposed for the solution of the minimum-fuel finite-thrust time-constrained rendezvous problem \cite{coverstone1994optimal, zavoli2014indirect}.
Indirect methods\cite{Bryson1979} are characterized by a high numerical accuracy and require a small computational effort.
However, the optimal mission structure, i.e., the sequence of burn and coast arcs, is usually unknown, and this poses a severe challenge on the capability of indirect methods to routinely, rapidly, and automatically solve the problem at hand\cite{ZavoliAlaska} (partially mitigated by the ongoing development of control regularization and continuation techniques)\cite{bertrand2002new}.
In a cooperative rendezvous, this issue is further stressed because the number of potentially-optimal mission structures grows drastically due to the simultaneous presence of two maneuvering spacecraft.
As a result, the definition of a suitable initial guess for the indirect method is hard and tedious, if not impractical at all.

On the other hand, direct methods, that rely on a transcription of the original problem into a general nonlinear programming problem, are easier to set up and are generally more robust to the initial guess. Nevertheless, the direct approach still requires a careful choice of the first guess because of the high-sensitivity to the tentative solution in complex problems. In addition, direct methods come with a greater computational burden and the optimality of the obtained solution is usually questionable.

Recently, convex programming has gained increasing popularity in the aerospace field thanks to the theoretical guarantees on the solution optimality and the availability of highly efficient numerical algorithms \cite{liu2017survey}. 
Even though most aerospace problems cannot be readily solved as convex optimization problems, several ideas have been proposed to convert a given nonconvex problem into a convex one, through a process referred to as convexification. 
Lossless convexification techniques exploit a change of variables in the problem formulation or a suitable constraint relaxation in order to obtain a convex problem having the same solution as the original problem \cite{accikmecse2011lossless}.
Successive convexification techniques, instead, rely on defining a sequence of convex subproblems  which leads to the original one.
The linearization of the equations of motion and of the nonconvex constraints around the previously found solution is the key of these techniques, which permit an efficient solution to realistic aerospace problems, usually involving nonlinear dynamics and/or nonconvex state and control constraints. 
The convergence properties of this approach have been proved in many applications~\cite{benedikter2019convexascent, lu2013autonomous, wang2016constrained}.

In complex problems undesired phenomena, such as artificial unboundness or artificial infeasibility, may compromise the solution procedure
and some safe-guarding expedients are necessary to ensure convergence.
In this regard, common approaches involve the use of a trust region or virtual controls in the convex formulation\cite{mao2016successive}.
The present paper proposes a novel filtering approach in the recursive update of the reference solution that aims at improving the algorithm robustness.
In particular, the reference solution, rather than simply being the last found, is computed as the weighted sum of the last \emph{three} solutions.
The principal merit of this technique is an easy implementation, especially when compared to other safe-guarding methods (e.g., adaptive trust regions). 
This filtering approach proved to be effective for the problem at hand, and its adoption in other optimization problem is worth to be investigated.

In the present paper a combination of lossless and successive convexification techniques is adopted for the solution of a cooperative time-fixed rendezvous problem.
Some nonlinearities are preserved in the convex subproblems by performing a convenient change of variables and a constraint relaxation, while a successive linearization of the equations of motion permits the use of the same dynamical model as the original (nonlinear) problem.
In order to show the effectiveness of the approach, numerical results for several cases are presented and compared with the solutions provided by an indirect method\cite{zavoli2014indirect}. 
In particular, the solver capability of finding the best solution in the presence of different families of optimal solutions, as detected by the indirect method, is discussed.

\section{Original Problem Formulation}

A finite-thrust time-constrained cooperative rendezvous is studied in this paper. 
At the initial time two identical spacecraft are moving on the same circular orbit of radius $\tilde{r}_0$ in phase opposition, i.e., with a phase difference equal to \num{180} degrees.
At the given final time $\tilde{t}_f$, they must meet on a circular orbit of assigned radius $\tilde{r}_f$ (greater than the initial one).
Both spacecraft can maneuver and share the same values for initial mass, maximum thrust, and specific impulse.
The goal of the optimization problem is to find the trajectories and the control laws that minimize the overall propellant consumption.

\subsection{System Dynamics}
In the present work, the spacecraft are modeled as point mass objects. Under this assumption, the state of each one is fully described by its position $\bm{r}$, velocity $\bm{v}$ and mass $m$. In order to favor an easier statement of the problem, the satellite position is given in polar Earth-Centered Intertial (ECI) coordinates, while the inertial velocity vector is expressed in a Local-Vertical-Local-Horizontal (LVLH) frame. Therefore, the state vector of each spacecraft is:
\begin{equation}
	\bm{x} = \begin{bmatrix}
	r & \theta & \varphi & v_r & v_t & v_n & m 
	\end{bmatrix}
	\label{eq:state_vector_original}
\end{equation}
where $r$ is the geocentric distance, $\theta$ is the right ascension, $\varphi$ is the latitude and $v_r$, $v_t$ and $v_n$ are, respectively, the radial, eastward and northward velocity components.
All variables are normalized with respect to the initial radius, the corresponding circular velocity and the mass of a single spacecraft.

In this study a simple inverse-square gravity model is considered. 
Indeed, the effects of the gravitational perturbations are deemed minimal. The gravitational acceleration is thus expressed as:
\begin{equation}
    \bm{g} = -\frac{\mu}{r^3} \bm{r}
    \label{eq:gravity}
\end{equation}
where $\mu$ is Earth's gravitational parameter.

The only control on each spacecraft is the thrust, that is expressed in the same frame as the velocity. So, the control variables for each satellite are:
\begin{equation}
    \bm{u} = \begin{bmatrix}
    T_r & T_t & T_n
    \end{bmatrix}
    \label{eq:thrust_vector_original}
\end{equation}
Both the thrust magnitude $T$ and direction $\bm{\hat{T}}$ must be determined by the optimization procedure. The thrust magnitude has an upper bound:
\begin{equation}
    T \leq T_{max}
    \label{eq:T_max}
\end{equation}
\begin{comment}
\hl{On the other hand}, $\bm{\hat{T}}$ must be a unit vector, hence the following condition must be ensured at any time:
\begin{equation}
	T_r^2 + T_t^2 + T_n^2 = 1
	\label{eq:thrust_unit_original}
\end{equation}
\end{comment}

In order to provide the thrust, the spacecraft must eject propellant mass over time.
Specifically, the mass flow rate is related to the engine thrust via the effective exhaust velocity $c = g_0 I_{sp}$, where $g_0$ is the gravity acceleration at sea level and $I_{sp}$ is the specific impulse in a vacuum:
\begin{equation}
    \dot{m} = -\frac{T}{c}
    \label{eq:mass_flow_rate}
\end{equation}

The resulting set of differential equations $\bm{\dot{x}} = \bm{f}(\bm{x}, \bm{u}, t)$ that describes the dynamics of a single spacecraft is the following:
\begin{align}
    \dot{r} &= v_r \label{eq:original_ODE_r} \\
    \dot{\theta} &= \frac{v_t}{r \cos\varphi} \label{eq:original_ODE_theta} \\
    \dot{\varphi} &= \frac{v_n}{r} \label{eq:original_ODE_varphi} \\
    \dot{v_r} &= \frac{v_t^2 + v_n^2}{r} -\frac{\mu}{r^2} + \frac{T_r}{m} \label{eq:original_ODE_v_r} \\
    \dot{v_t} &= -\frac{v_r v_t}{r} + \frac{v_t v_n}{r} \tan{\varphi} + \frac{T_t}{m} \label{eq:original_ODE_v_t} \\
    \dot{v_n} &= -\frac{v_r v_n}{r} - \frac{v_t^2}{r} \tan{\varphi} + \frac{T_n}{m} \label{eq:original_ODE_v_n} \\
    \dot{m} &= -\frac{T}{c} \label{eq:original_ODE_m}
\end{align}

\subsection{Objective Function}
The mission performance is evaluated in terms of the overall propellant consumption. 
Equivalently, since the initial mass of both satellites is assigned, the merit index to maximize can be defined as the sum of the spacecraft final masses:
\begin{equation}
    J = m_{\text{I}}(\tilde{t}_f) + m_{\text{II}}(\tilde{t}_f)
\end{equation}
where subscripts $\text{I}$ and $\text{II}$ refer to either satellite.

\subsection{Boundary Conditions}
In addition to the differential constraints, the optimization procedure has to take into account the boundary constraints.
The initial state of the two satellite is completely assigned, while the final state has to meet several mission requirements.

First of all, the rendezvous condition requires that both satellites share the same position and velocity at the final time. 
Some attention has to be paid while imposing this condition on the right ascension angle $\theta$. 
Indeed, this variable is not bounded in a limited interval (e.g., $\left[0, 2 \pi \right]$), thus it can assume any value in the range from zero to infinity.
Hence, while imposing the final rendezvous condition, one has to take into account an integer number of additional revolutions $k_{\text{rev}}$ performed by one of the satellites with respect to the other:
\begin{equation}
	\theta_{\text{I}}(\tilde{t}_f) = \theta_{\text{II}}(\tilde{t}_f) + 2 k_{\text{rev}} \pi \label{eq:rendezvous_condition_theta}
\end{equation}
For the specific case under investigation the optimal value of $k_{\text{rev}}$ is always zero.
So, there is no need to include it as an integer optimization variable.

Second, the two spacecraft must be on the same circular orbit of radius $\tilde{r}_f$ at the end of the mission. 
No other orbital element, but semi-major axis and eccentricity, is required for the target orbit. 
The final conditions are formulated as constraints on the position and velocity of either satellite:
\begin{align}
	r(\tilde{t}_f) &= \tilde{r}_f \label{eq:final_condition_r} 
\\
	v_r(\tilde{t}_f) &= 0 \label{eq:final_condition_v_r}
\\
	v_t^2(\tilde{t}_f) + v_n^2(\tilde{t}_f) &= \frac{\mu}{\tilde{r}_f} \label{eq:final_condition_v}
\end{align}
Actually, since the rendezvous conditions already ensures that the two satellites share final position and velocity, the conditions in Equations~(\ref{eq:final_condition_r}--\ref{eq:final_condition_v}) can be imposed only for one of the two spacecraft.

\section{Convex Transcription}
A convex optimization problem is characterized by a convex objective function, linear equality constraints, and inequality constraints that define a convex feasible set.
The problem under investigation cannot be readily solved by means of convex programming algorithms.
In fact, it has to be converted into a convex problem.
In this application, the original problem is converted into a special class of convex programming problems, the Second-Order Cone Programming (SOCP) problems.
A SOCP problem has a linear objective, linear equality constraints and second-order cone constraints. 
This class of programming problems allows for representing quite complex constraints and can be solved by means of highly-efficient interior point methods, even for a large number of variables\cite{alizadeh2003second}.

In this section, the convex optimization problem is formulated. 
First, lossless convexification is performed, then the remaining nonconvexities are handled by means of successive linearization.
Finally, the continuous-time convex problem is transcribed into a finite-dimensional problem.

\subsection{Change of Variables}
A change of variables is carried out to replace nonlinear terms in the dynamics by linear terms and obtain a control-affine dynamical system. 
This step is of crucial importance in order to favor a more robust successive convexification  process.
In addition, one should notice that the nature of the original problem is fully preserved, as no approximation is introduced.

The new control variables include both the thrust direction and the thrust-to-mass ratio. 
These are defined as:
\begin{align}
    u_r = \frac{T_r}{m} &&
    u_t = \frac{T_t}{m} &&
    u_n = \frac{T_n}{m} &&
    u_N = \frac{T}{m} \label{eq:new_controls}
\end{align}
By introducing $u_r$, $u_t$ and $u_n$ in Equations~(\ref{eq:original_ODE_v_r}--\ref{eq:original_ODE_v_n}), we obtain the following control-affine equations:
\begin{align}
    \dot{v_r} &= \frac{v_t^2 + v_n^2}{r} -\frac{\mu}{r^2} + u_r \label{eq:affine_ODE_v_r} \\
    \dot{v_t} &= -\frac{v_r v_t}{r} + \frac{v_t v_n}{r} \tan{\varphi} + u_t \label{eq:affine_ODE_v_t} \\
    \dot{v_n} &= -\frac{v_r v_n}{r} - \frac{v_t^2}{r} \tan{\varphi} + u_n \label{eq:affine_ODE_v_n}
\end{align}
However, introducing $u_N$ in Eq.~\eqref{eq:original_ODE_m} would not produce the same effect. A further step is required, and a new state variable is defined\cite{liu2018fuel}:
\begin{equation}
    z = \ln m
    \label{eq:z}
\end{equation}
Now, by differentiating Eq.~\eqref{eq:z} and combining it with Eq.~\eqref{eq:original_ODE_m} one obtains:
\begin{equation}
    \dot{z} = \frac{\dot{m}}{m} = -\frac{T}{m c} = - \frac{u_N}{c}
\end{equation}
that is an affine function of the control variables.

So, the new state and control vectors are:
\begin{gather}
    \bm{x} = \begin{bmatrix}
    r & \theta & \phi & v_r & v_t & v_n & z
    \end{bmatrix}
    \label{eq:state_vector_cvx} \\
    \bm{u} = \begin{bmatrix}
    u_r & u_t & u_n & u_N
    \end{bmatrix}
    \label{eq:control_vector_cvx} 
\end{gather}

\subsection{Constraint Relaxation}
The new defined control variables are related to each other by the following condition:
\begin{equation}
    u_r^2 + u_t^2 + u_n^2 = u_N^2
    \label{eq:thrust_direction_equality_path_con_new}
\end{equation}
that is a nonlinear, nonconvex, equality constraint.

A common convexification technique consists of relaxing a nonconvex constraint of the original problem into a convex constraint (\emph{constraint relaxation}).
Such an approach is particularly appealing since it does not introduce any approximation and allows to preserve the original problem nonlinearity in the convex formulation. 
This is a valuable feature since it has been shown that keeping some nonlinearities can significantly favor the convergence of a successive convexification algorithm\cite{yang2019comparison}.
In particular, the constraint relaxation is convenient especially when the optimal solution is the same as in the original problem, even though the relaxation defines a larger feasible set.
In such cases the relaxation is said to be \emph{exact}\cite{accikmecse2011lossless}.

The control constraint of Eq.~\eqref{eq:thrust_direction_equality_path_con_new} is suitable for an exact relaxation. So it can be safely replaced by the following inequality constraint:
\begin{equation}
	u_r^2 + u_t^2 + u_n^2 \leq u_N^2
	\label{eq:thrust_direction_cone_con}
\end{equation}
that corresponds to a second-order cone constraint. 

Finally, the newly defined $u_N$ variable should be bounded by the maximum thrust of the engine.
\begin{equation}
	0 \leq u_N \leq T_{max} e^{-z}
	\label{eq:u_N_path_con_nonconvex}
\end{equation}
However, this constraint is nonconvex and does not permit a convex relaxation, so it has to be linearized around a reference solution $(k)$ as:
\begin{equation}
	0 \leq u_N \leq T_{max} e^{-z^{(k)}} \left(1 - \left(z - z^{(k)}\right) \right)
	\label{eq:u_N_path_con_linearized}
\end{equation}

\begin{comment}
The constraints \eqref{eq:thrust_direction_cone_con} and \eqref{eq:u_N_path_con_linearized} form altogether a three-dimensional geometric half-cone in the control space, as shown in figure~\ref{fig:cone}.

\begin{figure}[!htb]
\centering
\includegraphics[height=60mm]{second_order_cone_u.pdf}
\caption{Relaxed control constraint results in a second-order cone constraint.}
\label{fig:cone}
\end{figure}

\end{comment}

The constraint of Eq.~\eqref{eq:u_N_path_con_linearized} permits the coexistence of burn,and coast arcs in a single phase. 
This is particularly advantageous for the problem under investigation since the burn structure of the solution is unknown \emph{a priori} and difficult to provide as satellites can maneuver at the same time.
The greatest drawback of mixing burn and coast arcs into a single phase is related to the fact that burn arcs require a much greater number of discretization nodes than coast ones in order to accurately approximate the continuous-time dynamics.
So, when using a uniform mesh, it has to be sufficiently dense everywhere, also in coast arcs, thus generating a large number of variables. 
This drawback is partially mitigated by the use of highly efficient numerical algorithms, that permit to solve even large problems with a limited computational effort.

\subsection{Successive Linearization}
At this point, the remaining nonconvexities are handled by successive linearization, i.e., the constraints are linearized around a reference solution $(k)$ that is iteratively updated.

By linearizing the equations of motion around a reference solution $(k)$ one obtains:
\begin{equation}
    \bm{\dot{x}} \approx \bm{f}\left( \bm{x}^{(k)}, \bm{u}^{(k)}, t \right) + \bm{f}_{\bm{x}}\left( \bm{x}^{(k)}, \bm{u}^{(k)}, t \right) \left( \bm{x} - \bm{x}^{(k)} \right) + \bm{f}_{\bm{u}}\left( \bm{x}^{(k)}, \bm{u}^{(k)}, t \right) \left( \bm{u} - \bm{u}^{(k)} \right)
    \label{eq:full_linearization}
\end{equation}
where $\bm{f}_{\bm{x}}$ and $\bm{f}_{\bm{u}}$ denote the partial derivative matrix of the right hand side of the equations of motion with respect to the state and control variables respectively.

Since the previous change of variables introduced a control-affine dynamics, the linearization reduces to:
\begin{equation}
    \bm{\dot{x}} = \bm{\hat{f}}\left( \bm{x}, t \right) + B\left( \bm{x}, t \right) \bm{u} \approx \bm{\hat{f}}\left( \bm{x}^{(k)}, t \right) + \bm{\hat{f}}_{\bm{x}}\left( \bm{x}^{(k)}, t \right) \left( \bm{x} - \bm{x}^{(k)} \right) + B\left( \bm{x}^{(k)}, t \right) \bm{u}
    \label{eq:affine_linearization}
\end{equation}
where $\bm{\hat{f}} = \bm{f}(\bm{u} = \bm{0})$ and $B$ is the $\left( n_x \times n_u \right)$ control coefficient matrix, with $n_x$ and $n_u$ denoting the number of state and control variables respectively.

By introducing the $\left( n_x \times n_x \right)$ $A$ matrix and the vector $\bm{c}$ of size $n_x$:
\begin{align}
    A &= \bm{\hat{f}}_{\bm{x}}\left( \bm{x}^{(k)}, t \right) \\
    \bm{c} &= \bm{\hat{f}}(\bm{x}^{(k)}, t) - A \bm{x}^{(k)} 
    \label{eq:c_vector}
\end{align}
the linearized dynamics becomes:
\begin{equation}
    \bm{\dot{x}} \approx A\left( \bm{x}^{(k)}, t \right) \bm{x} + B\left( \bm{x}^{(k)}, t \right) \bm{u} + \bm{c}\left( \bm{x}^{(k)}, t \right)
    \label{eq:matrix_linearization}
\end{equation}
where the coefficients of the dynamics depend only on the reference state $\bm{x}^{(k)}$ and not on the controls $\bm{u}^{(k)}$. This provides robustness to an iterative method, since intermediate controls, which may present high-frequency jitters \cite{liu2015entry}, do not affect the dynamics in the following iteration.

Finally, the boundary condition in Eq.~\eqref{eq:final_condition_v} has to be linearized. The result of the linearization is the following:
\begin{equation}
	{v_t^{(k)}}^2 + {v_n^{(k)}}^2 + 2 v_t^{(k)} (v_t - v_t^{(k)}) + 2 v_n^{(k)} (v_n - v_n^{(k)}) = \frac{\mu}{\tilde{r}_f} \label{eq:final_condition_v_lin}
\end{equation}

\subsection{Discretization}
\label{subsec:discretization}
The optimal control problem stated so far is infinite-dimensional since state and control variables are continuous-time functions. 
However, numerical methods for solving optimization problems require a \emph{finite} set of variables and constraints. 
% "For solving", verificato sul Betts
In order to convert the optimal control problem into a finite-dimensional problem, a direct transcription method is used.
% Si dice finite-dimensional, non finite-dimension (verificato sul Betts)

The independent variable, i.e., time, is discretized by dividing the mission duration into $M - 1$ intervals. So, one obtains $M$ points:
\begin{equation}
    t_0 = t_1 < \dots < t_{M} = \tilde{t}_f
    \label{eq:time_mesh}
\end{equation}
Each point $t_j$ is referred to as \emph{node} of the \emph{mesh}. 
Notice that, in general the nodes do not have to be equally spaced. 
In fact, in order not to generate a large discrete problem, the mesh should be dense only in the intervals where a small number of nodes would produce an inaccurate discrete approximation.

Once the grid is defined, both state and control variables are discretized over it. The differential constraints are replaced by a finite set of algebraic constraints, or \emph{defect} constraints. 
A simple trapezoidal integration scheme is employed in the current application, and the resulting defect constraint between nodes $j$ and $j + 1$ is:
\begin{equation}
    \bm{x}_{j + 1} - \left( \bm{x}_j + \frac{h_j}{2} \left(\bm{f}_j + \bm{f}_{j + 1}\right) \right) = \bm{0}
    \label{eq:trapz_ODE_con}
\end{equation}
Similarly, path constraints, such as Eq.~\eqref{eq:thrust_direction_cone_con}, are converted into a finite set of algebraic constraints by imposing them at each mesh node.

\subsection{Mesh Refinement}
The discrete-time problem is only an approximation of the original continuous-time problem.
Once the final solution is obtained, the quality of the discrete solution must be formally inspected and eventually a mesh refinement process has to be carried out in order to meet the desired tolerances.
Many techniques for the choice of the new mesh nodes have been proposed over the years.
The Betts and Huffman\cite{betts1998mesh} approach is used.

In order to evaluate the discretization error, the discrete solution must be converted into a continuous-time solution $(\bm{\Tilde{x}}(t), \bm{\Tilde{u}}(t))$ that approximates the real (unknown) solution $(\bm{\hat{x}}(t), \bm{\hat{u}}(t))$.
The state $\bm{x}(t)$ is approximated as a vector of cubic splines, with the conditions:
\begin{align}
    \bm{\Tilde{x}}(t_j) &= \bm{x}(t_j) \\
    \frac{d}{dt} \bm{\Tilde{x}}(t_j) &= \bm{f}\left( \bm{x}(t_j), \bm{u}(t_j), t_j \right)
\end{align}
Instead, the control is represented as a linear interpolation of the node values. 
Whereas the control is assumed to be correct and optimal, the error between the state $\bm{\Tilde{x}}(t)$ and the true solution is:
\begin{equation}
    \bm{\eta}_j = \int_{t_j}^{t_{j + 1}} |\bm{\Tilde{x}}(t) - \bm{\hat{x}}(t) | dt
    \label{eq:discretization_error}
\end{equation}
The integral in Eq.~\eqref{eq:discretization_error} can be estimated using a step size smaller than the one of the original grid. 
In particular, two trapezoidal (half) steps are used to estimate $\bm{\eta}_j$ as:
\begin{equation}
    \bm{\eta}_j \approx \frac{1}{2} \left| \bm{\Tilde{x}}(t_j + h_j) - \bm{\Tilde{x}}(t_{j}) - \frac{h_j}{4} \left( \bm{\Tilde{f}}_3 + 2 \bm{\Tilde{f}}_2 + \bm{\Tilde{f}}_1 \right) \right|
    \label{eq:discretization_error_estimate}
\end{equation}
where:
\begin{align}
\bm{\Tilde{f}}_k &= \bm{f}[\bm{\Tilde{x}}(s_k), \bm{\Tilde{u}}(s_k)] \\
s_k &= t_j + \frac{1}{2} (k - 1) h_j
\end{align}

If the error is above a given tolerance, the grid has to be refined. 
This is carried out by adding new nodes to the mesh. 
It has been observed that a basic refinement approach, such as simply taking twice as much intervals in each phase, may cause convergence problems of the successive convexification algorithm.
Hence, it is important to add as few points as possible. 
In Reference~\citenum{betts1998mesh} a method that selects new grid points by solving an integer programming problem is proposed. 
In particular, new points are selected to minimize the maximum discretization error by subdividing the current grid.
The initial reference solution on each new mesh is automatically obtained by interpolating the previous solution at the adjacent nodes.
This approach is deemed the most suitable for the problem under investigation as it guarantees that the grid size among successive refinements does not increase much, and only in the intervals above tolerance.

As a remark, the first mesh, on the one hand, must be sufficiently dense to approximate accurately the continuous-time problem and, on the other, it has to be coarse enough not to cause convergence problems\cite{kelly2015transcription}. Therefore, an adequate number of points must be picked while defining the starting grid.

\subsection{Successive Convexification Algorithm}

In order to converge towards the solution of the original problem a \emph{sequence} of SOCP problems has to be solved. 
Indeed, the SOCP problem formulated in the previous sections is only an approximation of the original problem as it considers the linearized dynamics rather than the real one.
Nevertheless, by updating iteratively the reference solution with the newly found solution, it has been shown in many applications that the process converges to the original problem solution.
The recursive process is terminated when the difference among the reference solutions goes below an assigned tolerance:
\begin{equation}
	\left\lVert \bm{x}^{(k)} - \bm{x}^{(k - 1)} \right\rVert_\infty < \epsilon_{\text{tol}}
	\label{eq:successive_cvx_termination_condition}
\end{equation}

Often, when employing successive linearization techniques, undesired phenomena, such as \emph{artificial unfeasibility} or \emph{artificial unboundedness}, may show up.
In these cases, safe-guarding modifications, such as virtual controls and a trust region, must be added to the convex formulation.
However, in the present application such phenomena did not show up, so neither virtual controls nor a trust region were necessary.

A different kind of phenomenon occurred when solving the problem under investigation. 
It has been observed that, even though usually convergence is attained in less than 10 iterations, in a few mission scenarios the reference solution tends to oscillate between two or more (non-acceptable) solutions, leading to an infinite process (bounded only by the limit on the maximum number of iterations).
In order to prevent such oscillations a \emph{filtering} technique is applied.

\begin{table}[htb]
    \caption{Reference solution update weights}
    \label{tab:weights_filtering}
    \centering
    \begin{tabular}{c c c}
        \hline
        $k_0$ & $k_1$ & $k_2$ \\
        \hline
        6/11 & 3/11 & 2/11 \\
        \hline
    \end{tabular}
\end{table}

The filtering concerns the recursive update of the reference solution, that, instead of being updated exclusively as the last obtained solution, is computed as a weighted sum of the \emph{three} previous solutions:
\begin{equation}
	x^{\text{ref}} = k_0 x^{(k)} + k_1 x^{(k - 1)} + k_2 x^{(k - 2)} 
\end{equation}
The values of the three weights used in the present application are reported in Table~\ref{tab:weights_filtering}. The proposed technique, extremely easy to implement, provides the required robustness to the successive convexification procedure.

\subsection{Initialization}
By definition, the local minimum of a convex problem is also its global minimum. Hence, contrarily to general nonlinear problems, the first guess choice does not affect the converged solution quality.
This is a unique feature that indirect and general NLP optimization methods do not share.
However, since successive convexification techniques were employed, the algorithm requires an initial reference solution.

In general, the reference solution for the direct optimization has to be as close as possible to the optimal solution, as the linearized constraints should accurately approximate the original ones.
The successive convexification, instead, in many applications can achieve convergence even starting from a trivial first guess.
For the problem under investigation a coasting trajectory on the initial orbit resulted to be effective.
Such reference is quite easy to provide as it requires no guess on the control variables nor on the burn arcs sequence of the two spacecrafts. 
So, actually, the initialization of the successive convexification algorithm does not represent an obstacle at all in the optimization process.

The easy initialization is one of the greatest advantages of the proposed approach.
Indeed, a more accurate reference trajectory would have been much more difficult to provide, as the contemporary presence of multiple maneuvering spacecraft greatly increases the complexity of the problem.
% that significantly favors the practical solution process, especially when compared with indirect and general NLP methods, that are highly sensitive with respect to the first guess.

\section{Coplanar Case}

In this and in the next section numerical results are presented in order to show the effectiveness of the proposed approach.
The presented algorithm has been implemented in C++ using Gurobi\cite{gurobi} as SOCP solver.
First, a planar (equatorial) case is investigated and the results are compared with those provided by an indirect method\cite{zavoli2014indirect}.
In the following section the study is extended to a non-coplanar case.

\begin{table}[htb]
    \caption{Spacecraft features}
    \label{tab:values}
    \centering
    \begin{tabular}{c c}
        \hline
        Quantity & Value \\
        \hline
        $T_{max}$ & \num{0.1} \\
        $c$ & \num{1} \\ 
        \hline
    \end{tabular}
\end{table}

At departure, the two satellites are on the same circular equatorial orbit of radius $\tilde{r}_0 = 1$ with a starting mass $m(t_0)$ equal to unit.
Since the two spacecraft are absolutely identical, satellite I is arbitrarily assumed to be the leading one ($\theta_{\text{I}}(t_0) = \pi$), and it is eventually reached by satellite II ($\theta_{\text{II}}(t_0) = 0$).
At the end time the two satellites must meet on a circular orbit of radius $\tilde{r}_f = 1.2$.
The spacecraft relevant features are reported in Table~\ref{tab:values}.

The mission requirements do not constrain the final orbit inclination.
However, since both satellites depart from the same orbital plane, the optimal target orbit  shares the same inclination.
Moreover, since the starting plane is equatorial, the latitude $\varphi$ and the normal velocity component $v_n$ should be null during the whole mission.
Therefore, the final condition of Eq.~\eqref{eq:final_condition_v} simply becomes:
\begin{equation}
	v_t(\tilde{t}_f) = \sqrt{\frac{\mu}{\tilde{r}_f}}
\end{equation}
This expression of the final velocity condition is particularly favorable as it can be readily included in the convex formulation without any linearization.

\begin{figure}[htb]
	\centering
	\includegraphics[width=0.7\linewidth]{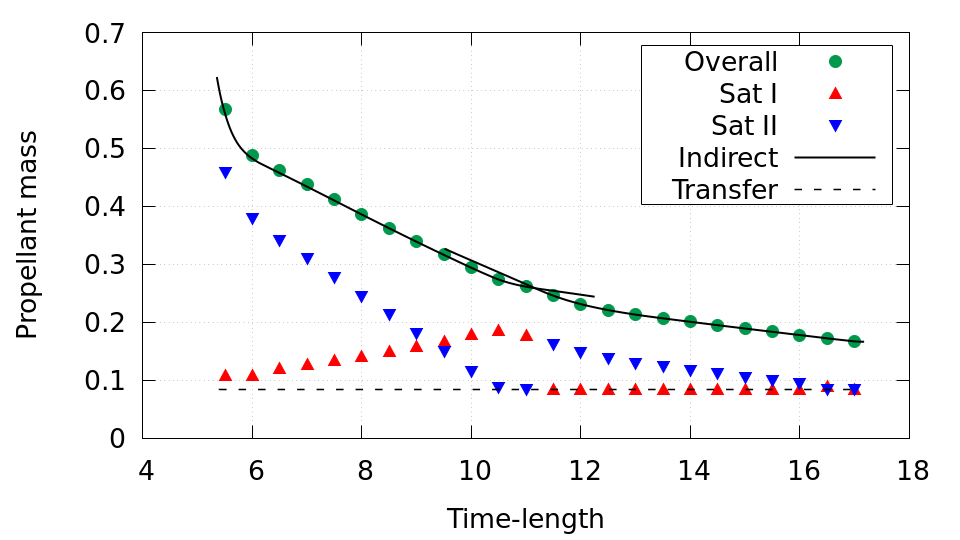}
	\caption{Propellant mass consumption versus mission duration for the coplanar case.}
	\label{fig:mp_Dt}
\end{figure}

The mission must be accomplished in a given time $\tilde{t}_f$.
If the available time is sufficiently long, both spacecraft perform a Hohmann-like maneuver, yet at different times in order to achieve the correct relative phasing.
Conversely, if the available time is not enough, the missions become more complex and expensive: only these missions are here considered.
Figure~\ref{fig:mp_Dt} reports the propellant consumption as a function of the available time.
The overall consumption provided by the presented approach (green symbols) is compared to analogous results (continuous line) obtained by an indirect method\cite{zavoli2014indirect}.
The indirect method detected two families of solutions that actually coexist in a limited interval of time-lengths.
As shown in Figure~\ref{fig:mp_Dt}, the convex approach converges towards the same solution as the indirect method for any mission time-length.
In particular, even when both families of solutions are viable, the convex algorithm successfully converges towards the best one, discarding automatically the most expensive one.

Figure~\ref{fig:mp_Dt} also reports the individual propellant contribution of each satellite.
It can be observed that in general the two spacecraft consume different amounts of propellant.
In particular, the roles of the two spacecraft are suddenly inverted at a time-length approximately equal to \num{11.165}, when a family of solutions becomes more convenient than the other.
From now on, the solutions for time-length smaller than \num{11.165} will be referred as belonging to family $\mathcal{A}$, while to the others as belonging to family $\mathcal{B}$.

\begin{figure}[htb]
\centering
\subfigure[$\tilde{t}_f = 10.5$]{\label{fig:familyA}\includegraphics[width=0.49\linewidth]{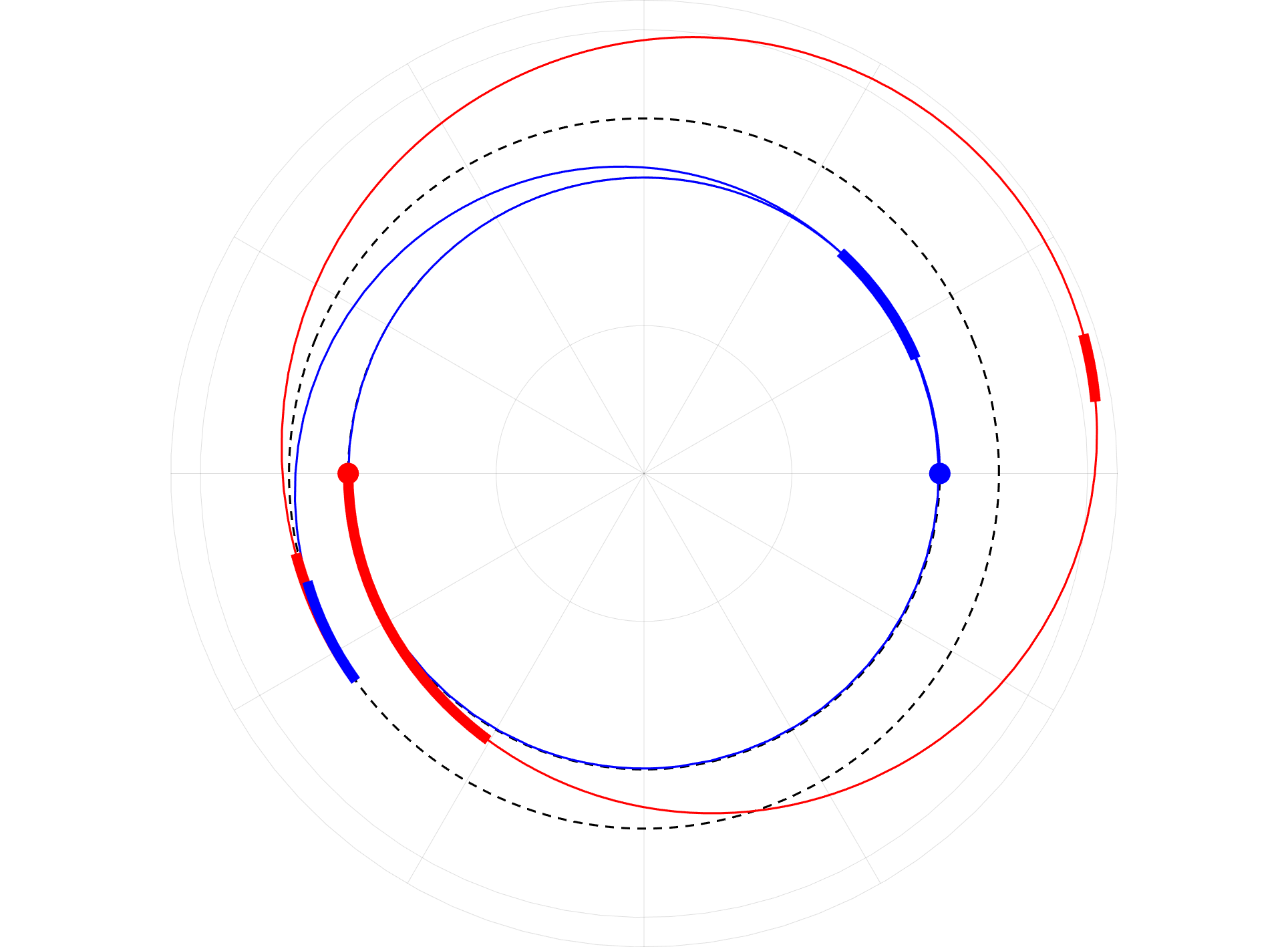}}
\subfigure[$\tilde{t}_f = 11.5$]{\label{fig:familyB}\includegraphics[width=0.49\linewidth]{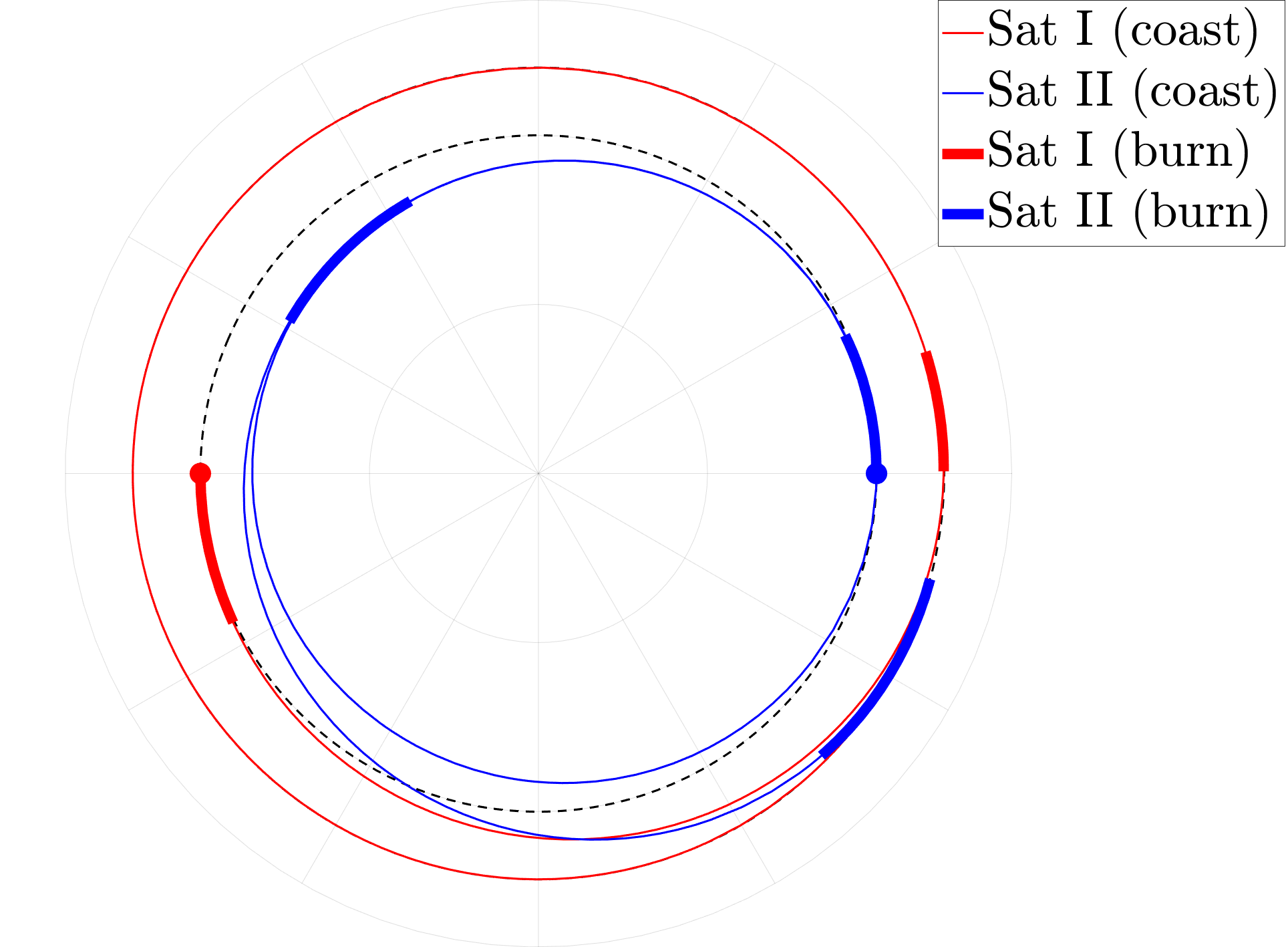}}
\caption{Trajectories of the two families (coplanar case)}
\label{fig:polar_traj_families}
\end{figure}

Figure~\ref{fig:polar_traj_families} presents the trajectories of the two families of solutions.
Family $\mathcal{A}$ (Figure~\ref{fig:familyA}) features the leading satellite (Sat I) performing approximately one revolution, while the other satellite performs one and a half revolutions.
Instead, family $\mathcal{B}$ solutions (Figure~\ref{fig:familyB}) are characterized by almost two complete revolutions of the trailing satellite. 
In the former case, the leading spacecraft flies over an external orbit while the other spacecraft stays on the initial orbit and recovers phase. 
In the latter case, the trailing spacecraft flies on a faster internal orbit and the other spacecraft can wait on the final orbit.

\subsection{Phasing duty}
The required mission effort can be divided into the cost to reach the final orbit (\emph{transfer cost} $\Delta m_T$) and the \emph{phasing duty} $\Delta m_\theta$, that is, the additional propellant consumed to achieve the rendezvous at the final time.
The transfer cost can be computed by optimizing the same problem without the $\theta$ rendezvous constraint of Eq.~\eqref{eq:rendezvous_condition_theta}.
Since all time-lengths are greater than the Hohmann transfer time, the transfer cost is the same for all the considered times, equal to \num{0.083} per satellite. 
Consequently, the phasing duty can be estimated as the difference between the propellant consumed in the full problem and the transfer cost:
\begin{equation}
	\Delta m_\theta = \Delta m - \Delta m_T \label{eq:}
\end{equation}

With reference to Figure~\ref{fig:mp_Dt}, where the dashed horizontal line represents the transfer cost of one satellite, for times longer than \num{17.17}, the phasing duty is null, and the mission cost is equal to the transfer cost of the two spacecraft only.
On the other hand, as the available time reduces, the overall phasing duty increases.
In solutions of family $\mathcal{B}$ the phasing duty is completely borne by satellite II. 
Indeed, satellite I does not have sufficient time to complete two revolutions on an orbit with a period greater than the final one;
hence, its trajectory is very close to the optimal Hohmann-like transfer followed by a coasting on the final orbit. 
Instead, spacecraft II starts its transfer by braking, then covers one and a half revolutions on an internal orbit, and at the second periapsis, accelerates to raise the apoapsis to the final circular orbit.

Solutions belonging to family $\mathcal{A}$ are characterized by a greater cooperation between the two satellites.
Indeed, the leading satellite raises its apoapsis over the desired one to wait for the trailing satellite, that in the meanwhile is flying on an internal orbit.
However, as the mission duration reduces, satellite I cannot help anymore satellite II because the available time permits a Hohmann-like transfer, but it is too short to allow for a further half revolution.
So, satellite II is forced to fly at lower altitudes, thus consuming a greater amount of propellant.

\subsection{Convergence Behavior}
For all the time-lengths investigated, the successive convexification algorithm achieves convergence.
In particular, already after the first 5 iterations, the intermediate solution closely resembles the final solution.
Indeed, the following iterations only refine the solution quality.
In most cases, the algorithm successfully terminates, i.e., the termination condition of Eq.~\eqref{eq:successive_cvx_termination_condition} is satisfied.
In a few other cases, some small oscillations are still detected, even though the solution filtering is active.
However, since the solution quality is acceptable already in the first iterations, an \emph{a priori} limit of 25 iterations resulted to be effective for the problem under investigation.

\begin{figure}[htb]
	\centering
	\includegraphics[width=0.7\linewidth]{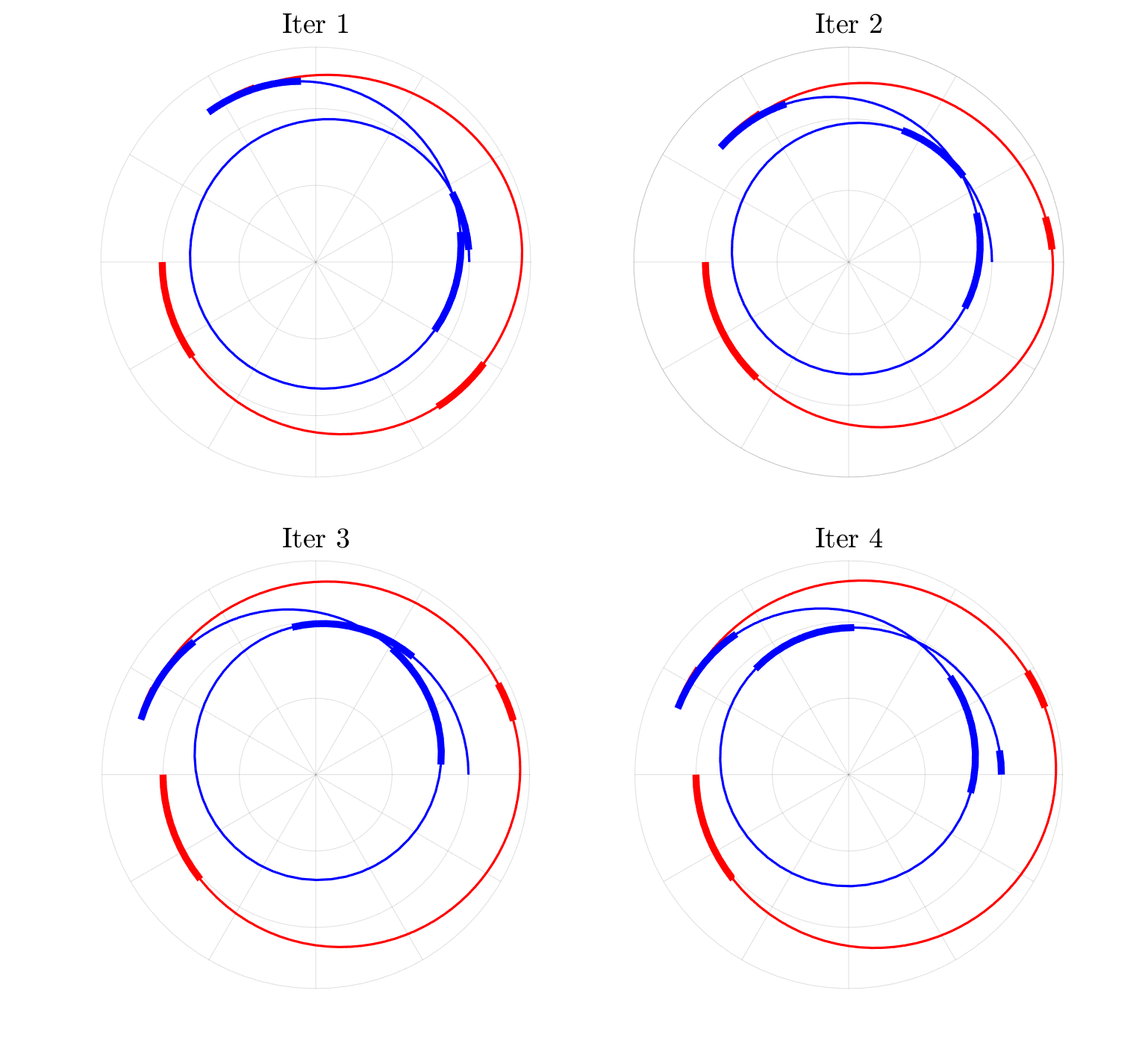}
	\caption{Spacecraft trajectories across a typical iteration sequence of the successive convexification algorithm}
	\label{fig:iterations}
\end{figure}

Figure~\ref{fig:iterations} shows an example of trajectories generated across a typical sequence of initial iterations.
One observes that the trajectory is basically defined at iteration 2 and remains almost constant in the next iterations.
Therefore, the state variables do not vary much across successive iterations.
Instead, the burn structure, thus the control variables, tends to change significantly among iterations.
This behavior would cause notable convergence problems to the successive convexification algorithm if the linearization depended upon the reference controls.
Instead, the dynamics was purposely transformed into a control-affine one, in order to prevent the convergence problems associated with control oscillations in intermediate iterations. 

Once convergence is attained, the quality of the discretization is inpected and eventually the mesh is refined.
The starting grid is made up of 101 equally spaced nodes.
This number of nodes permits to represent quite accurately the continuous-time problem and, at the same time, does not entail a large number of variables.
Furthermore, the chosen refinement algorithm ensures that new nodes are added only in intervals above tolerance, that mostly correspond to the burn arcs. 
This guarantees the problem to remain as small as possible, thus fast to solve also with the refined mesh.

As a final remark, it has been verified that the relaxation of the constraint of Eq.~\eqref{eq:thrust_direction_equality_path_con_new} into Eq.~\eqref{eq:thrust_direction_cone_con} is exact.
Indeed, even though the \emph{inequality} constraint is imposed, the control variables satisfy the equality constraint within tolerance.

\section{Non-coplanar Case}

%  demonstrating the trade-off between cost and time that characterizes the rendezvous problem. 

In this section the study is extended to a three-dimensional case.
The two satellites are assumed to depart from circular orbits of equal radius but different inclinations.
Satellite I is placed on an equatorial orbit, while satellite II is on an orbit with an inclination equal to 10 degrees.
At departure, both spacecraft are on the equatorial plane but at diametrically opposite points, that is, on the line of nodes. 
Satellite I is assumed to be the leading one ($\theta_{\text{I}}(t_0) = \pi$) while satellite II is the trailing one ($\theta_{\text{II}}(t_0) = 0$).

As in the previous case, no mission requirement on the final orbit inclination is imposed.
Now, since the two spacecraft depart from different orbital planes, the inclination of the arrival orbit is unpredictable and the final condition of Eq.~\eqref{eq:final_condition_v} cannot be simplified.

\begin{figure}[htb]
	\centering
	\includegraphics[width=0.7\linewidth]{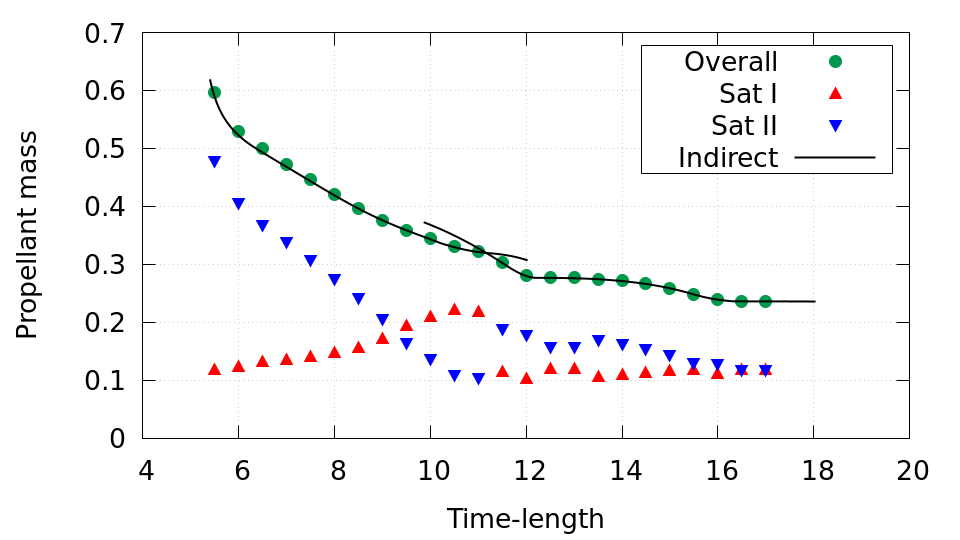}
	\caption{Propellant mass consumption versus mission duration for the non-coplanar case}
	\label{fig:mp_Dt_10deg}
\end{figure}

The mission time-length is assigned and if it is sufficiently long both spacecraft can perform a Hohmann-like maneuver.
However, now the maneuver must also include a plane change, that is quite an expensive maneuver.
In addition, the plane rotation can be executed only in certain points of the trajectories, thus it significantly affects the mission performance.
Figure~\ref{fig:mp_Dt_10deg} shows the propellant consumption as a function of the available time for the non-coplanar case.
The symbols used are the same as in Figure~\ref{fig:mp_Dt}.
The results of the convex approach are compared with those provided by an indirect method.
One can observe that the overall propellant consumption is practically the same in the two approaches.
Also this scenario features two families of solutions that coexist in a limited interval of time-lengths.
The convex approach successfully detects the most convenient family at every time.

\begin{figure}[htb]
\centering
\subfigure[$\tilde{t}_f = 10.5$]{\label{fig:familyA_3D}\includegraphics[width=0.45\linewidth]{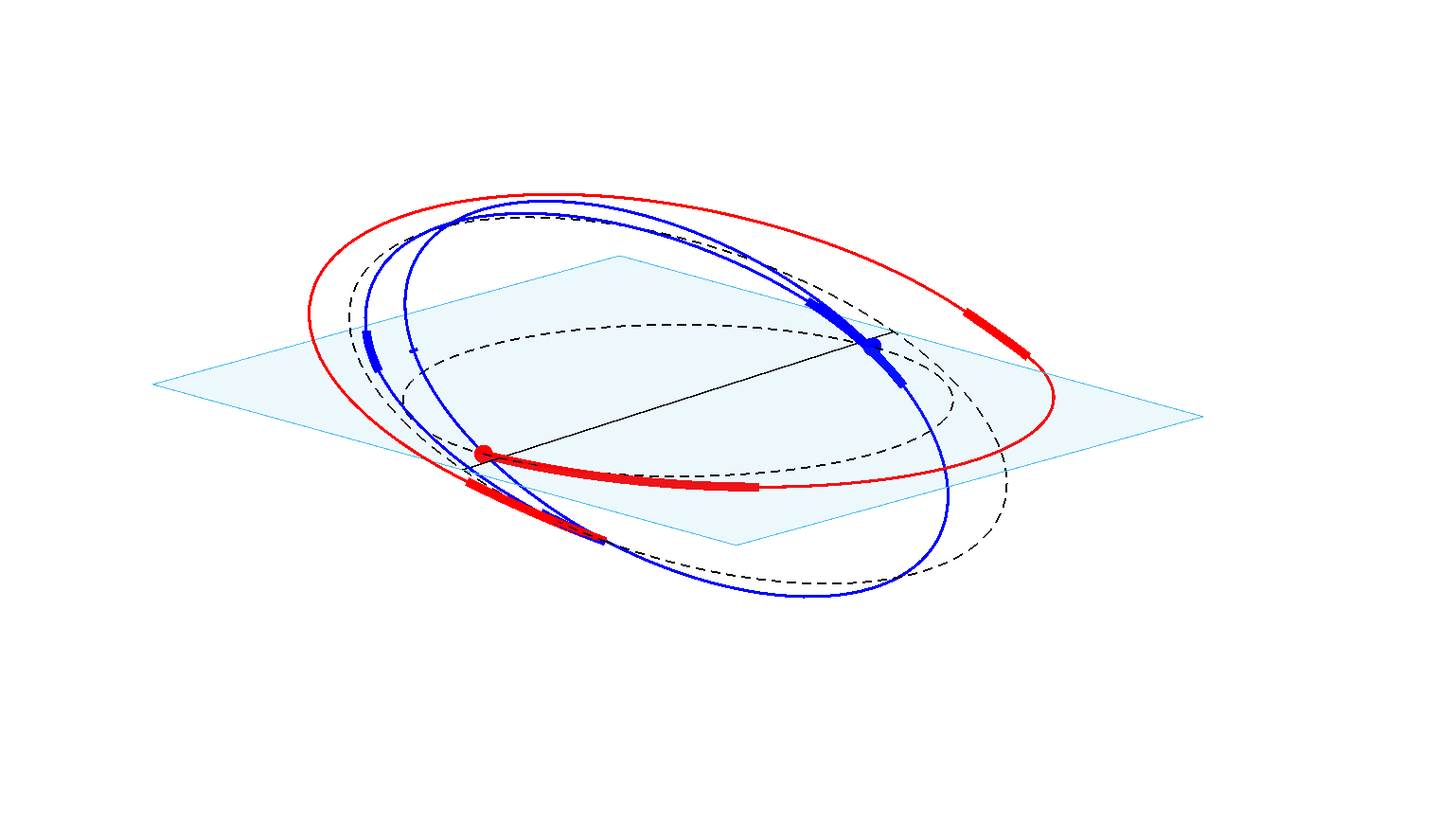}}
\subfigure[$\tilde{t}_f = 11.5$]{\label{fig:familyB_3D}\includegraphics[width=0.45\linewidth]{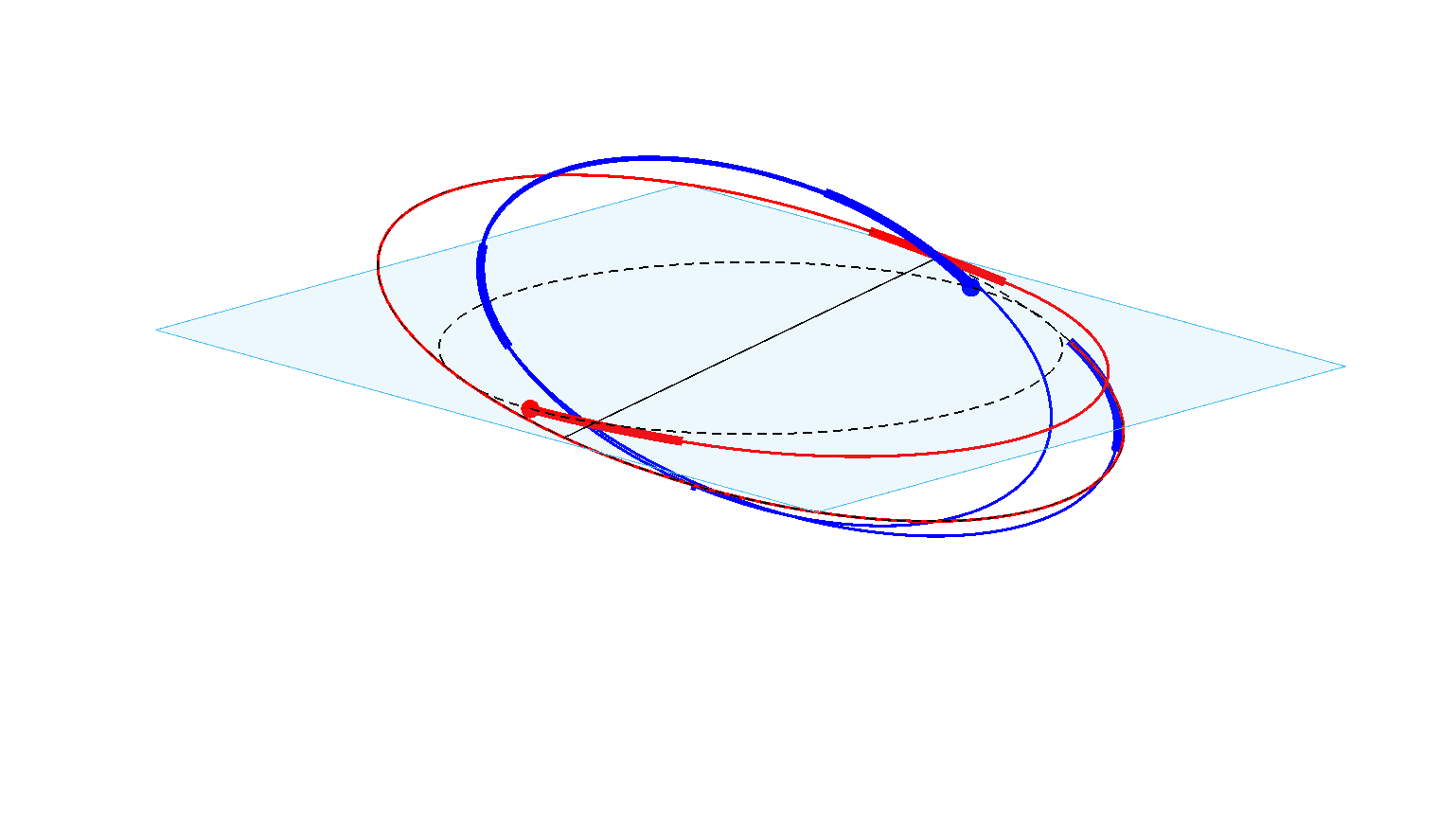}}
\caption{Trajectories of the two families (non-coplanar case)}
\label{fig:traj_families_3D}
\end{figure}

Figure~\ref{fig:traj_families_3D} illustrates the trajectories of the two families of solutions.
In order to ensure a good visibility the z-axis has been rescaled by a factor 3.
Apart from being three-dimensional, the two families of the non-coplanar problem closely resemble the ones of the two-dimensional problem.
Indeed, family $\mathcal{A}$ (Figure~\ref{fig:familyA_3D}), the one that is optimal for time-lengths shorter than \num{11.0}, still features only one revolution for the leading satellite on an external orbit and approximately one and a half revolutions for the trailing one on an internal orbit.
Instead, in solutions belonging to family $\mathcal{B}$ (Figure~\ref{fig:familyB_3D}) the trailing spacecraft performs almost two revolutions, whereas the leading satellite performs a Hohmann-like transfer and waits for the other spacecraft on the final orbit.

\begin{figure}[htb]
	\centering
	\includegraphics[width=0.7\linewidth]{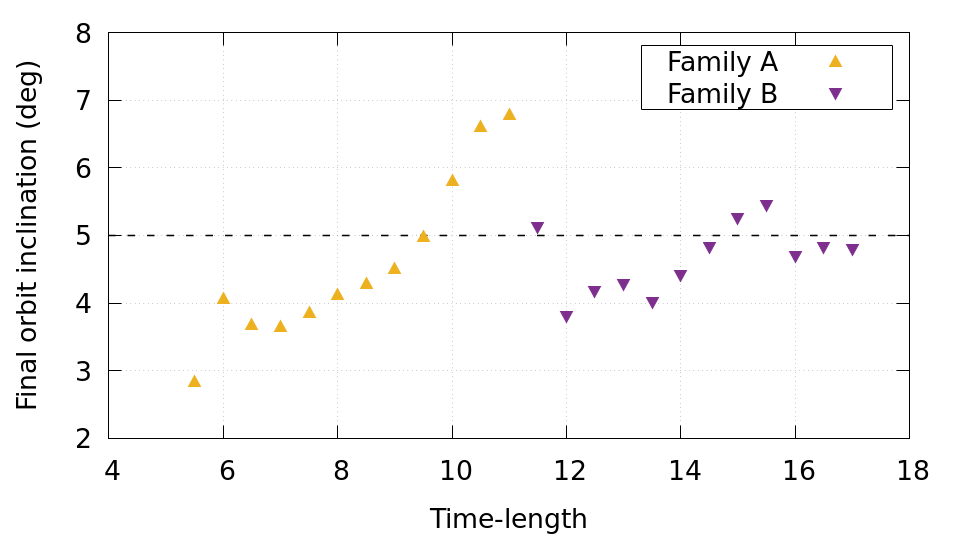}
	\caption{Final orbit inclination versus mission duration for the non-coplanar case}
	\label{fig:ainc_Dt}
\end{figure}

The main difference of the three-dimensional problem with respect to the coplanar one is the need to perform a plane change.
Whereas in the two-dimensional case the optimal target plane was equatorial, now the final plane depends on how the plane rotation maneuver is distributed between the two satellites.
The optimal final orbit inclination is reported in Figure~\ref{fig:ainc_Dt} as a function of the rendezvous time.
By comparing Figures~\ref{fig:mp_Dt_10deg} and \ref{fig:ainc_Dt} one can observe that in solutions of family $\mathcal{A}$ the satellite that maneuvers the most is also the one that performs the greatest plane change.
This repartition is optimal because of several reasons.
First, the spacecraft that executes most maneuvers is also the one that can combine in more occasions the rotation with the energy change.
Second, since it burns for longer periods, it is more likely to maneuver close to the nodes, where the plane rotation is more efficient.
Finally and most importantly, it quickly becomes the lightest satellite, hence it can consume less propellant while maneuvering.
Instead, the final inclinations of the solutions belonging to family $\mathcal{B}$ show some oscillations.
According to the considerations reported for family $\mathcal{A}$, one would have expected a more regular trend.
However, now the available time is longer, so the propellant consumption gets lower and another factor gains importance in the mission optimization: performing the plane rotation at the right moments.
Therefore, the inclination of the final orbit depends on which satellite can rotate the plane more efficiently, thus resulting in an unpredictable trend among different time-lengths.

% Within the solutions of family $\mathcal{A}$, the plane change cost gradually passes from being borne mostly by satellite II to satellite I. 
% Indeed, after $\tilde{t}_f = 9.0$ there is an inversion on the individual propellant contribution of each satellite, due to the final plane passing from inclination below the intermediate plane one to values above it.

Finally, the other considerations reported for the coplanar case extend also to the three-dimensional problem. 
In particular, the same convergence behavior has been observed and the constraint relaxation results to be exact also in this case.

\section{Conclusions}
This paper presented a convex methodology to solve the cooperative rendezvous problem.
Since the problem is not promptly suitable for a convex programming approach, the detailed convex transcription process has been described.
By using lossless and successive convexification techniques the real, nonlinear, dynamics of the problem is maintained.
Furthermore, a filtering technique, concerning the update of the reference solution, has been proposed to enhance the robustness of the successive convexification algorithm.
Such an expedient, albeit being very simple, revealed to be effective for the studied problem.

One of the principal merits of the proposed approach is the convergence towards the optimal solution even using a trivial first guess.
In fact, one of the major drawbacks of indirect methods is the need for a proper initial guess, that, for the problem under investigation, is quite difficult to provide as the mission scenario involves two maneuvering spacecraft, thus two switching control structures.
The convex approach has been proven capable of recognizing the most convenient solution also in the presence of multiple locally optimal families, as detected by the indirect method.
In addition, contrarily to general NLP methods, convergence is attained in a limited, short, time.

% Sviluppi futuri: orbite ellittiche, deployment cooperativo

Future work will include different requirements on the target orbit, for example, a highly elliptical orbit, as the convexification of the corresponding boundary conditions may be challenging.
Moreover, the same methodology can be applied to the similar problem of the cooperative deployment, that still involves multiple spacecraft that must achieve a desired relative phasing.
Finally, the proposed filtering method for updating the reference solution will be further studied and its effectiveness on other problems will be tested.

\bibliographystyle{AAS_publication}   % Number the references.
\bibliography{references}   % Use references.bib to resolve the labels.

\begin{thebibliography}{10}

\bibitem{coverstone1994optimal}
V.~Coverstone-Carroll and J.~E. Prussing, ``Optimal cooperative power-limited
  rendezvous between coplanar circular orbits,''  {\em Journal of Guidance,
  Control, and Dynamics}, Vol.~17, No.~5, 1994, pp.~1096--1102.

\bibitem{zavoli2014indirect}
A.~Zavoli and G.~Colasurdo, ``Indirect optimization of finite-thrust
  cooperative rendezvous,''  {\em Journal of Guidance, Control, and Dynamics},
  Vol.~38, No.~2, 2014, pp.~304--314.

\bibitem{Bryson1979}
A.~E. Bryson, Y.-C. Ho, and G.~M. Siouris, {\em Applied Optimal Control:
  Optimization, Estimation, and Control}, Vol.~9.
\newblock Institute of Electrical and Electronics Engineers ({IEEE}), 1979,
  10.1109/tsmc.1979.4310229.

\bibitem{ZavoliAlaska}
A.~Zavoli, F.~Simeoni, L.~Casalino, and G.~Colasurdo, ``Optimal cooperative
  deployment of a two-satellite formation into a highly elliptic orbit,''  {\em
  Advances in the Astronautical Sciences}, Vol.~142, 2012, pp.~3647--3663.

\bibitem{bertrand2002new}
R.~Bertrand and R.~Epenoy, ``New smoothing techniques for solving bang--bang
  optimal control problems—numerical results and statistical
  interpretation,''  {\em Optimal Control Applications and Methods}, Vol.~23,
  No.~4, 2002, pp.~171--197.

\bibitem{liu2017survey}
X.~Liu, P.~Lu, and B.~Pan, ``Survey of convex optimization for aerospace
  applications,''  {\em Astrodynamics}, Vol.~1, No.~1, 2017, pp.~23--40.

\bibitem{accikmecse2011lossless}
B.~A{\c{c}}{\i}kme{\c{s}}e and L.~Blackmore, ``Lossless convexification of a
  class of optimal control problems with non-convex control constraints,''
  {\em Automatica}, Vol.~47, No.~2, 2011, pp.~341--347.

\bibitem{benedikter2019convexascent}
B.~Benedikter, A.~Zavoli, and G.~Colasurdo, ``A Convex Approach to Rocket
  Ascent Trajectory Optimization,''  {\em 8th European Conference for
  Aeronautics and Space Sciences (EUCASS)}, 2019.

\bibitem{lu2013autonomous}
P.~Lu and X.~Liu, ``Autonomous trajectory planning for rendezvous and proximity
  operations by conic optimization,''  {\em Journal of Guidance, Control, and
  Dynamics}, Vol.~36, No.~2, 2013, pp.~375--389.

\bibitem{wang2016constrained}
Z.~Wang and M.~J. Grant, ``Constrained trajectory optimization for planetary
  entry via sequential convex programming,''  {\em AIAA Atmospheric Flight
  Mechanics Conference}, 2016, p.~3241.

\bibitem{mao2016successive}
Y.~Mao, M.~Szmuk, and B.~A{\c{c}}{\i}kme{\c{s}}e, ``Successive convexification
  of non-convex optimal control problems and its convergence properties,''
  {\em 2016 IEEE 55th Conference on Decision and Control (CDC)}, IEEE, 2016,
  pp.~3636--3641.

\bibitem{alizadeh2003second}
F.~Alizadeh and D.~Goldfarb, ``Second-order cone programming,''  {\em
  Mathematical programming}, Vol.~95, No.~1, 2003, pp.~3--51.

\bibitem{liu2018fuel}
X.~Liu, ``Fuel-optimal rocket landing with aerodynamic controls,''  {\em
  Journal of Guidance, Control, and Dynamics}, Vol.~42, No.~1, 2018,
  pp.~65--77.

\bibitem{yang2019comparison}
R.~Yang and X.~Liu, ``Comparison of Convex Optimization-Based Approaches to
  Solve Nonconvex Optimal Control Problems,''  {\em AIAA Scitech 2019 Forum},
  2019, p.~1666.

\bibitem{liu2015entry}
X.~Liu, Z.~Shen, and P.~Lu, ``Entry trajectory optimization by second-order
  cone programming,''  {\em Journal of Guidance, Control, and Dynamics},
  Vol.~39, No.~2, 2015, pp.~227--241.

\bibitem{betts1998mesh}
J.~T. Betts and W.~P. Huffman, ``Mesh refinement in direct transcription
  methods for optimal control,''  {\em Optimal Control Applications and
  Methods}, Vol.~19, No.~1, 1998, pp.~1--21.

\bibitem{kelly2015transcription}
M.~P. Kelly, ``Transcription methods for trajectory optimization,''  {\em
  Tutorial, Cornell University, Ithaca, New York}, 2015.

\bibitem{gurobi}
L.~Gurobi~Optimization, ``Gurobi Optimizer Reference Manual,''  2019.

\end{thebibliography}

\end{document}